\documentclass{article}

\usepackage{geometry}
\usepackage{lmodern} 
\usepackage{wrapfig}
\usepackage{colortbl}
\usepackage{amsmath}
\usepackage{amssymb}
\usepackage{mathrsfs}
\usepackage{makeidx}

\usepackage{multirow}
\usepackage{caption}
\usepackage{arydshln}

\newtheorem{thm}{Theorem}[section]

\usepackage{graphicx}

\usepackage{float}

\floatstyle{ruled}
\newfloat{triangle}{!ht}{tar}
\floatname{triangle}{arithmetic triangle}

\usepackage{ulem}

\usepackage{url}
\date{} 

\title{Generalized 3x + 1 Mappings : convergence and divergence}

\author{Robert Tremblay\\
				Boucherville,
				Canada (Qu\'ebec),\\				
				\texttt{roberttremblay02@videotron.ca}}

\begin{document}


\maketitle

\begin{abstract}
Discussion about the convergence and divergence of trajectories generated by certain functions derived from generalized $3x + 1$ mappings.
\end{abstract}

\section{Introduction}
In previous papers~\cite{tremblay_1, tremblay_2}, we analyzed the trajectories generated by the iterative application of several functions derived from generalized $3x + 1$ mappings. We have developed an algorithm that allow us to determine the necessary conditions for the existence or not of loops (cycles). From the periodicity property associated with the different trajectories for a given length, it has been possible to demonstrate that the number of cycles is finite. This led us to the notion of convergence or not of the trajectories towards these cycles which can be either closed or open. The function that gives rise to the original Collatz problem produces nine closed cycles and, since the number of cycles is limited, all the integers not belonging to these cycles are in infinite trajectories. The function that gives rise to the $3x + 1$ problem produces four opened cycles with the negative integers and the zero, and only one appears for positive integers. All other natural numbers seem to converge towards this cycle. The function generating the $5x + 1$ problem seems to lead to trajectories convergent and divergent. 

In this paper we will analyze the convergence and the divergence of trajectories associated with the $3x + 1$ and $5x + 1$ problems. The approach we use here is independent of the results previously obtained if we accept the conjecture which states that the number of cycles is finite (\cite{lagarias},\cite{matthews}).

We will use an intrinsic property, called periodicity, resulting from the iterative application of the functions. \textsl{A priori}, the behavior of trajectories generated for these functions seems chaotic. In fact, if we make appropriate groupings of trajectories we quickly observe a regularity in their distribution. In this context, we will be able to provide answers to the questions of convergence or not of these trajectories.


In this new version of this paper, we have added a section in the problem $3x+1$ called "Distribution function $F(k)$". We start with the approach undertaken by Riho Terras~\cite{terras} and taken up by several authors including Lagarias~\cite{lagarias} which allows to develop a density function $F(k)$. Terras proves that this function is well defined and has very interesting properties. We use his reasoning until the remarkable result which appears with the theorem of periodicity. Thereafter we use a completely different path which includes the properties generated by only two theorems developed in our paper. Our reasoning is much simpler and leads to results which coincide with those expected. In our opinion, this last exercise reinforces the power of the use of diophantine equations for the solution of problems like the two treated in this paper.    
        

\section{Mathematical stools}  

Mappings can be define on integers represented by functions such that each element of the set $\mathbb{Z}$ is connected to a single element of this set. The iterative application of these functions produces a sequence of integers called trajectories. In this paper, we will use indifferently the term \textsl{trajectory} or simply \textsl{sequence} to indicate a sequence of integers generated by these functions. We can easily construct the equations connecting any two integers of the sequence. The result can always be expressed in the general form $c = ax + by$, where $x$ represents any integer and $y$ is the integer resulting from the iterative application of the function. The parameters $a$ and $b$ depend on the function itself. So, we have diophantine equations of first degree at two unknowns. From a well-known result of this theory, we have the theorem

\begin{thm}
Let the diophantine equation $c = ax + by$ of first degree at two unknowns. If the coefficients $a$ and $b$ of $x$ and $y$ are prime to one another (if they have no divisor other than $1$ and $-1$ in common), this equation admits a infinity of solutions to integer values. If $(x_0, y_0)$ is a specific solution, the general solution will be $(x = x_0 + bq, y = y_0 - aq)$, where $q$ is any integer, positive, negative or zero.
	\label{equationDiophantine}
\end{thm}

\textsl{Proof}

References : on the web and~\cite{bordelles} 

$\blacksquare$ 

The distribution of integers resulting from the different groupings of trajectories will follow different progressions, and several of them will be of geometric type. Write a general geometric series as 

\begin{equation}
	\sum_{k=1}^{n}ar^{k-1} = ar^0 + ar^1 + ar^2 + ar^3 + \cdots + ar^{n-1},
\label{geometricSeries}	
\end{equation}

where $ar^{n-1}$ is the nth term of the series.

The sum of first $n$ terms is given by

\begin{equation}
	S_n = \frac{a(1-r^n)}{1-r}.
\label{geometricSum}	
\end{equation}

The terms of the geometric series will be represented by multinomials. 

We have the binomial formula

\begin{equation}
	\left( x + y \right)^n = \sum_{k=0}^{n} \left( \begin{array}{c} n \\ k \end{array}\right) x^ky^{n-k},
\label{binomial}	
\end{equation}

where  

\begin{equation}
	BC_{n,k} = \left( \begin{array}{c} n \\ k \end{array}\right) = \frac{n!}{k!(n-k)!},
\label{binomial_coefficients}	
\end{equation}

are called binomial coefficients $BC_{n,k}$.

The multinomial formula (Wikipedia) is

\begin{equation}
	\left( x_1 + x_2 + \cdots + x_m \right)^n = \sum_{k_1 + k_2 + \cdots + k_m}^{n} \left( \begin{array}{c} n \\ k_1, k_2, \cdots, k_m \end{array}\right) \prod_{t=1}^{m} x_t^{k_t},
\label{multinomial}	
\end{equation}

where  

\begin{equation}
	MC = \left( \begin{array}{c} n \\ k_1, k_2, \cdots, k_m \end{array}\right) = \frac{n!}{k_1!k_2! \cdots k_m},
\label{multinomial_coefficients}	
\end{equation}

is a multinomial coefficient. The sum is taken over all combinations of nonnegative integers indices $k_1$ through $k_m$ such that the sum of all $k_i$ is $n$. That is, for each term in the expansion, the exponents of the $x_i$ must add up to $n$.

For example, the third power of the trinomial $a + b + c$ is given by

\begin{equation*}
	\left( a + b + c\right)^3 = a^3 + b^3 + c^3 + 3a^2b + 3a^2c + 3b^2a + 3b^2c + 3c^2a + 3c^2b + 6abc,
\label{example_multinomial_coefficients}	
\end{equation*}

where

\begin{equation*}
	a^2b^0c^1 \phantom{1} has \phantom{1} the \phantom{1} coefficient \phantom{1} \left( \begin{array}{c} 3 \\ 2, 0, 1 \end{array}\right) = \frac{3!}{2!\cdot 0!\cdot 1!} = 3,
\label{multinomial_coefficient_1}	
\end{equation*}

and, 

\begin{equation*}
	a^1b^1c^1 \phantom{1} has \phantom{1} the \phantom{1} coefficient \phantom{1} \left( \begin{array}{c} 3 \\ 1, 1, 1 \end{array}\right) = \frac{3!}{1!\cdot 1!\cdot 1!} = 6.
\label{multinomial_coefficient_2}	
\end{equation*}

\section{Problem 3x + 1}
 
In first, we present three functions which encode the $3x + 1$ problem. 

Let the \textsl{Collatz function} $C(n)$ be defined as follow  

\begin{equation}
C(n)=\left\lbrace  
\begin{array}{ll}
3n+1 & \mbox{, if $n\equiv1\pmod{2}$}\\
\\
\frac{n}{2} & \mbox{, if $n\equiv0\pmod{2}$}.\\
\end{array}
\right.  
\label{Collatz_function}
\end{equation}

Even though this function, as well as the next two, is valid for all integers $n$, positive, negative or zero, we will use the set of natural numbers in most of the examples that follow. 

In a $2012$ paper~\cite{delahaye}, Delahaye gives a very good introduction to the $3x + 1$ problem. He produced a figure, that he called a tree, giving the directional tracking of several positive integers leading to the number $1$. We see that each integer has $1$ antecedent or $2$, never more. Some odd integers such as $3$, $9$, $15$, $\cdots$, are not preceded by any odd integers.  Several representations of trees are found on the web, in particular that giving the trajectories less than $20$ before reaching the number $1$~\cite{wikipedia_1}. In fact, the root of these trees is $1$, and they are built from inverse algorithms to those generated by the Collatz function. If the conjecture that states that all natural numbers end on the cycle $\langle1, 4, 2\rangle$ is true, these trees must cover all natural numbers.

Kontorovich and Lagarias~\cite{kontorovich_lagarias} in a paper produced in $2009$ work with two other functions with iterations faster than the Collatz function.

The first is the \textsl{$3x + 1$ function} $T(n)$ (or $3x + 1$ map)

\begin{equation}
T(n)=\left\lbrace  
\begin{array}{ll}
\frac{3n+1}{2} & \mbox{, if $n\equiv1\pmod{2}$}\\
\\
\frac{n}{2} & \mbox{, if $n\equiv0\pmod{2}$}.\\
\end{array}
\right.  
\label{3x+1_function}
\end{equation}

This function results from the fact that the $3x + 1$ operation applied to any odd integer always give an even integer. 

The second function, the \textsl{accelerated $3x + 1$ function} $U(n)$, is defined on the domain of all odd integers, and removes all powers of $2$ at each step. It is given by 

\begin{equation}
	U(n) = \frac{3n + 1}{2^{ord_2(3n+1)}}
\label{accelerated_function}
\end{equation}

in which $ord_2(n)$ counts the number of powers of $2$ dividing $n$. The function $U(n)$ was studied by Crandall in $1978$~\cite{crandall}. 

          
\subsection{Groupings (in triplets) of trajectories generated by the function $U(n)$}

Let us represent the \textsl{accelerated $3x + 1$ function} $U(n)$ by the following $3$ operations

\begin{equation}
U(n)=\left\lbrace  
\begin{array}{lll}
\frac{3n+1}{4} & \mbox{, if $n\equiv1\pmod{8}$}\\
\\
\frac{3n+1}{2} & \mbox{, if $n\equiv3\pmod{4}$}\\
\\
\frac{3n+1}{8} & \mbox{, if $n\equiv5\pmod{8}$}\\
\\
\end{array}
\right.  
\label{accelerated_function_2}
\end{equation}

The first two operations result in odd integers and the last operation ends with an alternation of even and odd integers $2\pmod{3}$. The table~\ref{Triplet_function_U} we give the three groupings resulting of these operations including the intermediate (abreviate 'int') operation $(3x + 1 )/2$. In each grouping we have a sequence of three integers ('triplets'). The intermediate integer is repeated in the second grouping. The reason will appear on its own. 

Even if the results contained in the table are easily deduced during its construction, it is possible to find them by writing the diophantine equations and, using the theorem~\ref{equationDiophantine}.

For example, the first operation leads to the equation

\begin{equation*}
	y = \frac{3x + 1}{4} \phantom{1234} or \phantom{1234} 1 = 4y - 3x.	
\end{equation*}

The couple of values $(x_0 = 1; \phantom{1} y_0 = 1)$ is a solution. According to the theorem \ref{equationDiophantine}, the general solution will be $(x = 1 + 4q_1; \phantom{1} y = 1 + 3q_1)$, where $q_1$ is any integer, positive, negative or zero. If this first operation result in odd integers, we can write $y = 2q_2 -1$. Then,

\begin{equation*}
	y = 1 + 3q_1 = 2q_2 -1 \phantom{1234} and \phantom{1234} q_1 = \frac{2q_2 - 2}{3}.	
\end{equation*}

The solutions to the integer values are $q_2 = 1, 4, 7, 10, \cdots$ and $q_1 = 0, 2, 4, 6, \cdots$.

We will have, as expected, 

\begin{equation*}
x = 1 + 4q_1 = 1, 9, 17, 25, \cdots \phantom{1234}  or \phantom{1234} x\equiv1\pmod{8}, 
\end{equation*}  

and

\begin{equation*}
y = 1 + 3q_1 = 1, 7, 13, 19, \cdots \phantom{1234}  or \phantom{1234} y\equiv1\pmod{6}. 
\end{equation*} 

The same procedure is used to find the solutions of the other operations.







If we analyze the decomposition of trajectories in triplets, it is possible to perceive a certain regularity. We will not further develop the search for the regularity by this way. We will instead define an another function, so $\mathfrak{U}(n)$, that apply to all integers. This new function produces trajectories whose triplet decomposition will be such that each of the integers of these triplets are in correspondence one-to-one with each of the elements of one and only one triplet produced by the function $U(n)$. We can then build a tree that, \textsl{a priori}, will include all the integers and bring out the regularity.

          
\subsection{New function $\mathfrak{U}(n)$ and one-to-one correspondence with the function $U(n)$}

We define the function $\mathfrak{U}(n)$

$$\mathfrak{U}(n)=\left\lbrace
\begin{array}{ll}
\frac{3n+1}{4} & \mbox{, if $n\equiv1\pmod{4}$}\\
\\
\frac{3n}{2} & \mbox{, if $n\equiv2\pmod{2}$}\\
\\
\frac{n+1}{4} & \mbox{, if $n\equiv3\pmod{4}$}.\\
\end{array}
\right.$$

This function does not come from generalized $3x + 1$ mappings. The table~\ref{Triplet_function_U_transformed} give three groupings resulting of these operations. In fact, this table brings together the trajectories of all integers and the integer resulting from a single application of the function $\mathfrak{U}(n)$ (duos). We form triplets by repeating the first integer. We can easily verify that all the integers $n$ of triplets in the table~\ref{Triplet_function_U} (function $U(n)$) are calculated using two simple operations on all the integers $n_{(new)}$, where $n_{(new)}$ is any integer beginning the triplets in the table~\ref{Triplet_function_U_transformed} (function $\mathfrak{U}(n)$), 

\begin{equation}
	n = 2n_{(new)} - 1 
\label{first_operation}
\end{equation} 

and 

\begin{equation}
	n = 3n_{(new)} - 1. 
\label{second_operation}
\end{equation}

All intermediate integers $2\pmod{3}$ of triplets of $U(n)$ and the last of the third grouping come from the operation $(3n_{(new)} - 1)$. All other integers (odd integers) come from the operation $(2n_{(new)} - 1)$. Then, the correspondence between the results of the application of the function $\mathfrak{U}(n)$ on all integers is one-to-one with the application of the function $U(n)$ on all odd integers.

As we will see quite rapidly, this new function will show the regularity in the distribution of sequences of integers. It will allow us to build a tree in which all the branches fit into each other and, possibly containing all integers. 

          
\subsection{Tree generated by the function $\mathfrak{U}(n)$ and distribution of integers}

Let us define the branches as the sequence of the integers that are connected to each other by the operations $(3n + 1)/4$ or $3n/2$ starting with $2\pmod{3}$ and ending with an integer resulting from the operation $(n + 1)/4$. The branches are interrelated, thus forming a tree (table~\ref{Tree}) whose the trunk is represented by the sequence $2 \rightarrow 3 \rightarrow 1$. If there is only one cycle ($1 \rightarrow 1$) for the natural numbers, then all the integers of a branch are different two by two, and also are two by two different from all the integers of any other branch, except of course, the end-of-branch integer that connect to another branch. 

We have retained the first and the last integer of the branch of length $19$ starting with $14$ and ending with $56$, at the top of the tree. For the branch starting with $56$ and whose length is $13$, we selected 4 values, namely the first ($56$) and the last three.

The search for the convergence of all integers towards $1$ is to determine if all the integers are in branches (distribution) interlocking into each other and if, sooner or later, they converge to smaller integers.

We will prove in a first time that all integers are found in sequences with an end of all possible lengths and their distribution is characterized by a geometric progression. 

Using the operation $(3n_{(new)} - 1)$ on all integers in branch ends of the tree in the table~\ref{Tree} and the operation $(2n_{(new)} - 1)$ on all the other integers of this tree, we build the tree represented in the table~\ref{Tree_odds} generated by the function $U(n)$. Nevertheless, it is easier to analyze the function $\mathfrak{U}(n)$ which applies to all integers and which naturally reveals the regularity in the distributions of the branches. 

          
\subsubsection{Distribution of integers in sequences with an end}

The solution of diophantine equations used below is easily verified from the groupings of table~\ref{Triplet_function_U_transformed}. The sequences that we will analyze end with an integer resulting the operation $(n + 1)/4$, but the beginning of the sequences is any integer, not necessarily a branch beginner integer $2\pmod{3}$.   

Let's first look the sequences of two integers whose second is obtained from the first by the operation $(n + 1)/4$ (end of branch). The length of the sequences is $L = 2$. The diophantine equation generating these sequences is $y = (x + 1)/4$, or $1 = 4y - x$ (general form is $c = by + ax$), of which a particular solution is $(x_0 = 3, y_0 = 1)$. By the theorem~\ref{equationDiophantine} the general solution will be $(x = 3 + 4q, y = 1 + 1q)$. $b = 4$ and $-a = 1$ are the increments which are added to $x_0$ and $y_0$ to get the general solution. In notation with modulo we will write $x\equiv3\pmod{4}$ and $y\equiv1\pmod{1}$. This result tells us that at every $4$ consecutive integers there is $1$ which ends after the operation $(n + 1)/4$ and this is repeated periodically to infinity. The distribution is $1$ for $4$, so $1/4$. Each integer $x$ in this group of length $2$ goes to a smaller integer. In resume, for the sequences of length $L = 2$ we build a single diophantine equation and the distribution of the first integers $x$ of these sequences is $1/d$ with $d = 4$ the denominator in the first expression of the diophantine equation.  

Let the sequences of three integers (length $L = 3$). Two cases are possible, whether the operations $3n/2$ (first type) or $(3n + 1)/4$ (second type) is applied to an integer followed by the operation $(n + 1)/4$. The resolution of two diophantine equations leads to two solutions, $x\equiv2\pmod{8}$ ($y\equiv1\pmod{3}$) and $x\equiv9\pmod{16}$ ($y\equiv2\pmod{3}$). The distributions are $1/8$ and $1/16$, meaning respectively that at each $8$ consecutive integers $x$ there is one that starts a sequence of a first type, and at each $16$ consecutive integers there is one that starts a sequence of a second type and this is repeated periodically to infinity. 

We can do this for $L = 4, 5, 6 , \cdots$. 

The results below are valid  for $L \ge 2$.

The total number of diophantine equations $\#(DE)_L$  is

\begin{equation}
	\#(DE)_L = 2^{L-2}. 
\label{number_DE}
\end{equation}

The different increments (or modulo) for the integers $x$ starting these sequences are

\begin{equation}
	Increment = d_1*d_2*d_3*\cdots*d_{L-1},
\label{increment}
\end{equation}

where $d_1, d_2, d_3, \cdots, d_{L-2} = 2,4$ and $d_{L-1} = 4$. 

The increment (or modulo) for the last integers $y$ is $3^{L-2}$.

The different distributions of first integers $x$ are calculated from the following binomial

\begin{equation}
	D_L = \frac{1}{4}\left(\frac{1}{2} + \frac{1}{4}\right)^{L-2},
\label{geometricSeries_firstDistribution}	
\end{equation}

which is the nth term ($n = L-1$) of a geometric progression. From equation~\ref{geometricSeries}, which is the general form for this series, we put $a = \frac{1}{4}$ and $r = \frac{3}{4}$. 

The binomial coefficients $BC_{n,k}$ of this binomial give the number of different diophantine equations for a given increment.

The sum of $n$ first terms is given by the equation~\ref{geometricSum}

\begin{equation}
	S_L = \frac{1}{4}\left(\frac{(1-(\frac{3}{4})^{L-1})}{(1-\frac{3}{4})}\right) = 1-\left(\frac{3}{4}\right)^{L-1}.
\label{geometricSeries_firstDistribution_sum}	
\end{equation} 

When $L$ goes to infinity the sum tends towards one, meaning that {\bf all integers $x$, without exception, start sequences with end}. 

In table~\ref{Sequences_L} and the table~\ref{Distribution} we give the sequences and distributions of first integers $x$ generated by the function $\mathfrak{U}(n)$ for the lengths $L = 2$ until $L = 6$ (and until $L = 7$ in the table~\ref{Distribution}). 

According to the equation~\ref{geometricSeries_firstDistribution_sum}, a large part of natural integers is included in the first sequences. The sequences with $L = 2$ until $L = 20$ include almost $99.6 \%$ of all natural integers.

The next step in our approach is to determine the part of natural integers that go to smaller integers.

The first sequences where the first integers $x$ are smaller than the final integers $y$ ($x<y$) are those starting with $x\equiv48\pmod{64}$ for the lengths $L = 6$. In fact, these sequences are the only ones of all those for the lengths $L = 2$ until $L = 6$ with this behavior. All other integers for sequences covering these lengths end with smaller integers ($x>y$). The distribution $1/64$ meaning that at each $64$ consecutive integers there is one that starts a sequence which end with a larger integer. 

For the sequences of length $L = 7$ and more we observe the same type of behavior. Some sequences with the smaller modulo will have the beginning integer smaller than that the integer of the end of the sequence. For all others we have the inverse behavior. It was quite predictable.

\begin{thm}
Let the sequences (of length $L$) of the integers that are connected to each other by the operations $(3n + 1)/4$ or $3n/2$ and ending with an integer resulting from the operation $(n + 1)/4$. The diophantine equation connecting the first integer $x$ and the last integer $y$ of  a sequence can be expressed in the general form $c = by - ax$ where the parameters $a$, $b$ and $c$, always positive, depend on the operations themselves and in which orders they are applied. If $b > a$, $x \ge y$ and, if $b < a$, $x < y$. 
	\label{Distribution_PP_PG}
\end{thm}

\textsl{Proof}

Let $k_1, k_2 = 0, 1, 2, \ldots$ and $k_3 = k_1 + k_2 = L - 2$, with $L \ge 2$.

Then, $a = 3^{k_3}$, $b = 4 \cdot 2^{k_1} \cdot 4^{k_2}$ and $c > 0$.

As the factors $a$ and $b$ of $x$ and $y$ are prime to one another, the diophantine equation admits a infinity of solutions to integer values. If $(x_0, y_0)$ is a specific solution, the general solution will be $(x = x_0 + bq, y = y_0 + aq)$, where $q$ is any integer, positive, negative or zero.

Two cases are possible, $b > a$ or $b < a$.

\vspace{2mm} 

\underline{First case : $b > a$}

Suppose that a particular solution $(x_0, y_0)$ is such that $x_0 < y_0$. We have the general solution

\begin{equation*}
	y = y_0 + aq \phantom{1234} and \phantom{1234} x = x_0 + bq,	
\end{equation*}

where $q$ is any integer, positive, negative or zero. As $b > a$ and $x_0 < y_0$, beyond a certain value of $q$, we will have $x > y$. The equation

\begin{equation*}
	c = by - ax,	
\end{equation*}

eventually lead to a negative $c$ value. But, the parameter $c$ must always be positive. Therefore $x > y$ when $b > a$.

\vspace{2mm}

For example, for $L = 2$, we have $k_1 + k_2 = L -2 = 0$, $k_1 = k_2 = 0$ and $k_3 = k_1 + k_2 = 0$. Then, $a = 3^{k_3} = 3^0 = 1$ and $b= 4 \cdot 2^{k_1} \cdot 4^{k_2} = 4$. We write the diophantine equation

\begin{equation*}
	y = \frac{x + 1}{4} \phantom{1234} or \phantom{1234} 1 = 4y - x.	
\end{equation*}

where $b = 4$, $a = 1$ and $c = 1$.

The couple of values $(x_0 = 3; \phantom{1} y_0 = 1)$ is a particular solution. The general solution will be $(x = 3 + 4q; \phantom{1} y = 1 + q)$. 

\vspace{2mm}

\underline{Second case : $b < a$}
 
 Let the equation

\begin{equation*}
	c = by - ax.	
\end{equation*}

As $c$ is always positive and $b < a$, $x$ must necessarily always be smaller than $y$ ($x < y$).

\vspace{2mm}

The sequences of length $L = 6$ beginning with $x\equiv48\pmod{64}$ treated in the following subsection are typical examples of this case. We have $b = 64$, $a = 81$ and $c = 16$. A particular solution is $(x_0 = 48, y_0 = 61)$, and the general solution is $(x = 48 + 64 q, y = 61 + 81q)$. Therefore $b < a$ and $x < y$. The table~\ref{Distribution_x_y} gives the values of $b$ and $a$ for $L = 2$ until $L = 8$ (function $\mathfrak{U}(n)$).

$\blacksquare$   

In table~\ref{Distribution_smaller} we give the first distributions of the natural numbers whose first of the sequence is smaller than the last ($x<y$). For sequences lengths $L = 2$ until $L = 20$ there is around $12 \%$ of {\bf all integers} (actually $99.58 \%$) behaving this way. On the other hand, around $88 \%$ of all natural numbers have their first numbers larger than the last of the sequence ($x>y$). This first result is quite remarkable because it tells us that not only $88 \%$ of all natural numbers go to smaller integers, but this is done in very specific slices (increments). 

We have all the elements necessary to build the algorithm of a function $\mathfrak{F}$ determining the fraction $f$ (or the \%) of natural integers such as $x < y$. For $L = 2$ to $20, 30, 40, 50, 60$,  

\begin{equation*}
	f = 0.1198839, 0.1236245, 0.1238451, 0.1238577, 0.1238584.  
\end{equation*}

If we take consecutive integer slices of $5 \phantom{1} 000$, for example, we can already see that the number of integers whose first of the sequence is smaller than the last is close of $12 \%$ (using an appropriate algorithm). The results (by slices) are indeed those anticipated by the function $\mathfrak{F}$.


A good part of the natural integers seem to go towards smaller integers. We will reinforce this fact by analyzing a little more in detail the some $12 \%$ of integers with $x < y$. 

          
\subsubsection{Distribution of integers in sequences with an end at the second level and more}

Let is now the sequences whose beginning integer $x$ is smaller than the end integer $y$, that is $12 \%$ of all integers. The first sequences with this behavior are those of length $6$, so $x\equiv48\pmod{64}$. This sequences are carried out by applying 4 times the operation $3n/2$ and ending with the operation $(n+1)/4$. The diophantine equation is written as

\begin{equation*}
	\left(\left(\frac{3}{2}\right)^4x+1\right)\cdot\frac{1}{4} = y \phantom{1234} or \phantom{1234} 16 = -81x + 64y.
\label{example_1}	
\end{equation*}

Let $(x_0 = 48, y_0 = 61)$ be a particular solution. Then, $(x = 48 + 64q, y = 61 + 81q)$ is the general solution, where q is any integer, positive, negative or zero. If we continue the iteration of the function $\mathfrak{U}(n)$ on the last integer of the sequence until reaching a second end of branch, we will have for example,

\begin{flushleft}
(48$\rightarrow$72$\rightarrow$108$\rightarrow$162$\rightarrow$243$\rightarrow$61)$\rightarrow$(61$\rightarrow$46$\rightarrow$69$\rightarrow$52$\rightarrow$78$\rightarrow$117$\rightarrow$88$\rightarrow$132$\rightarrow$198$\rightarrow$297$\rightarrow$223$\rightarrow$56)\quad PP

(112$\rightarrow$168$\rightarrow$252$\rightarrow$378$\rightarrow$567$\rightarrow$142)$\rightarrow$(142$\rightarrow$213$\rightarrow$160$\rightarrow$240$\rightarrow$360$\rightarrow$540$\rightarrow$810$\rightarrow$1215$\rightarrow$304)\quad PP

(176$\rightarrow$264$\rightarrow$396$\rightarrow$594$\rightarrow$891$\rightarrow$223)$\rightarrow$(223$\rightarrow$56)\quad PG

(240$\rightarrow$360$\rightarrow$540$\rightarrow$810$\rightarrow$1215$\rightarrow$304)$\rightarrow$(304$\rightarrow$456$\rightarrow$684$\rightarrow$1026$\rightarrow$1539$\rightarrow$385)\quad PP

\ldots	
\end{flushleft} 

where $PP$ ("Plus Petit") indicates that the beginning integer $x_1$ of the first sequence is smaller than the last integer $y_2$ of the second sequence, and $PG$ means larger ("Plus Grand"). 

The first sequence will be called the first level sequence and the second, the second level sequence. Then, $x_1\equiv48\pmod{64}$ and $y_1\equiv61\pmod{81}$. $y_1$ becomes then the beginning integer $x_2$ of the second sequence.

In table~\ref{Second_level_sequences} we have the first second level sequences $x_2\equiv61\pmod{81}$. We find a distribution identical to that obtained for the sequences of first level. All integers $61\pmod{81}$ start sequences with an end, of all possible lengths and their distribution is characterized by a geometric progression. 

For example, let's take sequences of length $L = 2$. The diophantine equation is

\begin{equation*}
	\frac{(x_2 + 1)}{4} = y_2,
\label{example_2a}	
\end{equation*}

where $x_2 = 61 + 81q$. Then,

\begin{equation*}
	62 = -81q + 4y_2.
\label{example_2b}	
\end{equation*}

A specific solution at integer values is $(q_0 = 2, y_{2_0} = 56)$ and the general solution is $(q = 2 + 4k,y_2 = 56 + 81k)$. We replace $q = 2 + 4k$  in $x_2 = 61 + 81q$. Then,

\begin{equation*}
	x_2 = 223 + 324k,
\label{example_2c}	
\end{equation*}

as expected. 

We proceed in the same way for all possible trajectories starting with $x_2 = 61 + 81q$. In writing the diophantine equations $c = ax + by$, we will have the coefficient $a = -81$ and the coefficient $b$ resulting from the products of $2$ and $4$. According to theorem~\ref{equationDiophantine}, since $a$ and $b$ are prime between them, the equations always have an infinity of solutions to the integers values. All trajectories starting with $x_2\equiv61\pmod{81}$ are possible and their distribution will be such that nearly $88 \%$ of them will see the beginning integer $x_2$  larger than the last integer $y_2$. 

An identical result would have been obtained for each of the $12 \%$ of sequences where $x_1< y_1$. Then, around $88 \%$ of all these sequences have their first integers $y_1 = x_2$  larger than $y_2$. This does not mean that $88 \%$ of $x_1$ is larger than $y_2$. For example, the sequence starting with $x_1 = 48$ end with $y_1 = 61$ at the first level and $y_2 = 56$ at the second level. At the second level, $x_2 = y_1 = 61$ and $y_2 = 56$ ($x_2 > y_2$, but $x_1 < y_2$). In fact, around $60 \%$ of $12 \%$ of $x_1$ is larger than $y_2$. 

As the sequences of integers at the second level develop similarly to those of the first level, we can use the function $\mathfrak{F}$ that we noted at the end of the previous section and call it recursively. For sequences up to lengths $30$ at the first level (for $L1 = 2$ to $30$), we had

\begin{equation}
	f = 0.1236245.  
\label{fraction_level_1}	
\end{equation}

At the second level (for $L1, L2 = 2$ to $30$) $f$ becomes  

\begin{equation}
	f = 0.05112079.  
\label{fraction_level_2}	
\end{equation}

With the third level (for $L1, L2, L3 = 2$ to $30$),

\begin{equation}
	f = 0.024040812.  
\label{fraction_level_3}	
\end{equation}

We have programmed the function $\mathfrak{F}$ in VBA (Visual Basic for Applications) with a standard laptop (standard computer). With the level $4$, the operations preformed by calling the function increase very rapidly and the integers appearing there become very large (over 29 digits). By using a more powerful computer, we could obtain the results of the function for levels higher than $3$.       

With the second level, we now have around $95 \%$ of all natural numbers that end up on a smaller integer and, $97.6 \%$ with the third level. This percentage increases with the fourth level and more. Not only $97.6 \%$, or more, of all natural numbers go to smaller, but this is done in very specific slices. 

In table~\ref{Distribution_by_slices} we give the results by slices ($1$ to $5 \phantom{1} 000, 10 \phantom{1} 000,  100 \phantom{1} 000, 1 \phantom{1} 000 \phantom{1} 000$), so the fraction of integers with $x < y$ for the first $7$ levels, beginning with the first natural number $1$. We can verify that the results would have been similar if we had chosen slices starting with any integer other than $1$. In the last column we find the first three values of $f$ obtained from the programmed function $\mathfrak{F}$ (equations~\ref{fraction_level_1}, \ref{fraction_level_2}, \ref{fraction_level_3}). Once again, the results (by slices) are indeed those anticipated by the function $\mathfrak{F}$.

          
\subsection{Distribution of odd integers in sequences with an end (function $U(n)$)}

We use a method similar to the one previously used to find the distribution of odd integers in the trajectories generated by the function $U(n)$. 


In table~\ref{Sequences_L_odds} we give the sequences of first odd integers $x$ generated by the function $U(n)$ for the lengths $L =2$ to $L = 6$. These results can be found by writing the different diophantine equations and solving them. Nevertheless, and because the two functions ($U(n)$ and $\mathfrak{U}(n)$) are one-to-one correspondence, we can simply use the integers of the trajectories constituting the branches given in the previous table~\ref{Sequences_L}, and apply the transformation $3n - 1$ to the end integers of the sequences and the transformation $2n - 1$ to all other integers.

The first trajectories for which the first odd integer $x$ is smaller than the last integer $y$ ($x < y$) is $x\equiv15\pmod{64}$ and $y\equiv20\pmod{81}$ for $L = 5$, corresponding to $x\equiv8\pmod{32}$ and $y\equiv7\pmod{27}$ for the function $\mathfrak{U}(n)$ where $x$ is larger than $y$ ($x > y$).

The complete distribution of all odd integers at the first level will be slightly different from the distribution of all integers with the function $\mathfrak{U}(n)$ in the table~\ref{Distribution_smaller}. Nearly $20 \%$ (instead of $12 \%$) of {\bf all odd integers $x$, without exception}, start sequences with end, such as $x < y$. Then, $80 \%$ of all odd integers go to smaller integers at the first level. 

At the second level, we have about $60 \%$ of $20 \%$ going to smaller integers (whose $x > y$). With the second level, we now have around $92 \%$ of all odd integers that end up on a smaller integer and this percentage increases with the third level and more. Not only $92 \%$, or more, of all odd integers go to smaller, but this is done in very specific slices.

In table~\ref{Distribution_by_slices_U} we give the results by slices, so the fraction of integers with $x < y$, beginning with the first odd integer $1$. In the last column we find the first three values of $f$ obtained from the programmed function $\mathfrak{F}$.
   
          
\subsection{Distribution function $F(k)$}

Let us define the distribution function $F(k)$ as

\begin{equation}
	F(k) = \lim\limits_{m \rightarrow \infty} \left( 1/m \right) \mu \{ n \le m \phantom{1} | \phantom{1} \chi(n) \ge k \},
\label{Function_distribution}	
\end{equation} 

where $\mu$ is the number of positive integers $n \le m$ with $m$ that tends towards infinity. $\chi(n)$ is called the "stopping time", and corresponds to the smallest positive integer such that the iterative application ($k$ times) of $3x+1$ function $T$ (equation~\ref{3x+1_function}) on a integer $n$ gives the result $T^kn<n$.

Terras~\cite{terras} proves that this function is well defined for any value of $k$ and that it tends towards $0$ for $k$ tending towards infinity.

Lagarias~\cite{lagarias} redoes the demonstration using the function we will call $G(k)$,

\begin{equation}
	G(k) = \lim\limits_{x \rightarrow \infty} \frac{1}{x} \phantom{1} \# \{ n : n \le x \phantom{1} and \phantom{1} \sigma(n) \le k \},
\label{Function_distribution_Lagarias}	
\end{equation} 

where $\sigma(n)$ is the "stopping time". This function $G(k)$ is in away almost the reciprocal of the function $F(k)$, and tends towards $1$ when $k$ tends towards infinity. The properties inherent in these functions will be clarified in the following examples.

Terras~\cite{terras} and Everett~\cite{everett} have independently demonstrated an important theorem (periodicity) bringing out the regularity in the trajectories generated by the $3x+1$ function $T(n)$. We extended this property to the original Collatz problem~\cite{tremblay_2} and the correspondent function $g(n)$. To do this, we used a fundamental property of the diophantine equations instead of the induction used by the cited authors.    
 

The periodicity theorem can be interpreted as follows.

Let $k$ be a number of iterations applied to any $2^k$ consecutive integers. We will have all possible combinations $2^k$ of operations $n/2$ on the even integers and $(3n+1)/2$ on the odd integers of the diadic sequences generated by the function $T(n)$ and each combination appears only once. All the integers $m$ of the form $m = n + 2^k$ will have the same combination of operations. The distribution of different combinations is then binomial versus the operations.

For example, let $k=1$ and the $2^k=2^1=2$ consecutive positive integers $1$ and $2$. The sequences of length $k + 1 = 2$ generated by the function $T(n)$ will be

\begin{equation*}
	(1,2) \phantom{1} (2,1) \phantom{1} (3,5) \phantom{1} (4,2) \phantom{1} (5,8) \phantom{1} (6,3)  \phantom{1} \cdots,
\label{sequences_2_numbers}	
\end{equation*} 

where we have added the numbers $3$, $4$, $5$ and $6$ after the two consecutive numbers so as to bring out the periodicity.

If we use the diadic sequences of the $0$ and $1$ representing respectively the even and odd operations, we will have

\begin{equation*}
	(1) \phantom{1} (0) \phantom{1} (1) \phantom{1} (0) \phantom{1} (1) \phantom{1} (0) \phantom{1} \cdots,
\label{sequences_2_numbers_diadic}	
\end{equation*} 

all repeating periodically for every two consecutive sequences. This result follows from the fact that the all integers alternate between the even and odd integers. 

Let's write the first diophantine equation for the even integers which give the first integer $x$ of the sequence versus the last integer (here the second) of the sequence,

\begin{equation*}
	\frac{x}{2} = y ;\phantom{1234} 0=2y-x ;\phantom{1234} and \phantom{1} the \phantom{1} general \phantom{1} solution \phantom{1234} (2 + 2q, 1 + 1q).	
\end{equation*}

The second equation is the same as the first, but in the form $c=by-ax$ with $c=0$, $b=2$ and $a=1$. Since the parameters $b$ and $a$ are prime to each other, the theorem~\ref{equationDiophantine} allows us to write the general solution $(x_0 +b q, y_0 + aq)$, where $(x_0, y_0)$ is a particular solution and $q$ is any integer, positive, negative or zero. According to the theorem~\ref{Distribution_PP_PG}, as $b > a$ then $x \ge y$. In this case, the stopping time is equal to the number of iterations $k = 1$.  

The second diophantine equation for the odd integers is,

\begin{equation*}
	\frac{3x+1}{2} = y ;\phantom{1234} 1=2y-3x ;\phantom{1234} and \phantom{1} the \phantom{1} general \phantom{1} solution \phantom{1234} (1 + 2q, 2 + 3q).	
\end{equation*}

Here, $c=1$, $b=2$ and $a=3$. According to the theorem~\ref{Distribution_PP_PG}, as $b < a$ then $x < y$. In this case, the stopping time is greater than the number of iterations $k = 1$. Then, all sequences starting with an odd positive integer contribute to the distribution function $F(k)$ because $\chi > k$. As expected, $F(k=1) = 1/2$.


Let another example. Take $k=2$ and the $2^k=2^2=4$ consecutive positive integers $3,4,5$ and $6$. The sequences of length $k + 1 = 3$ generated by the function $T(n)$ will be

\begin{equation*}
	(3,5,8) \phantom{1} (4,2,1) \phantom{1} (5,8,4) \phantom{1} (6,3,5) \phantom{1} (7,11,17) \phantom{1} (8,4,2)  \phantom{1} \cdots,
\label{sequences_4_numbers}	
\end{equation*} 

where we have added the numbers $7$ and $8$ after the four consecutive numbers so as to bring out the periodicity. 

The diadic sequences are

\begin{equation*}
	(1,1) \phantom{1} (0,0) \phantom{1} (1,0) \phantom{1} (0,1) \phantom{1} (1,1) \phantom{1} (0,0) \phantom{1} \cdots,
\label{sequences_4_numbers_diadic}	
\end{equation*} 

all repeating periodically for every four consecutive sequences. 

We can write the $2^k = 2^2 = 4$ diophantine equations in the same way as before. But, we will do it differently here. In fact the diadic sequences we will help to deduce whether or not the stopping time is greater or less than the number of iterations $k = 2$.

In the general case, the parameter $b = 2^k$ and the parameter $a = 3^{k_2} \cdot 1^{k_1} = 3^{k_2}$ with $k$ the total number of iterations, $k_1$ the number of operations on the even integers, and $k_2 = k - k_1$ the number of operations on the odd integers. 


\begin{table}[H]
\begin{center}
\begin{tabular}{|c|c|c|c|c|c|c|c|}
	
	\hline
	$diadic$ & $k_1$ & $k_2$ & $b = 2^k$ & $b = 3^{k_2}$ & $b \phantom{1} vs \phantom{1} a$ & $x \phantom{1} vs \phantom{1}  y$ & $stopping$ \\
	
	$sequences$ & & & & & & & $time \phantom{1} \chi(n)$ \\
	\hline
	$(0,0)$ & 2 & 0 & 4 & 1 & $b > a$ & $x > y$ & $-$ \\
	
	$(0,1)$ & 1 & 1 & 4 & 3 & $b > a$ & $x > y$ & $-$ \\
	
	$(1,0)$ & 1 & 1 & 4 & 3 & $b > a$ & $x > y$ & $\chi = k$ \\
	
	$(1,1)$ & 0 & 2 & 4 & 9 & $b < a$ & $x < y$ & $\chi > k$ \\
	\hline
\end{tabular}
\end{center}
\caption{Stopping time for k = 2}
\label{stoppingTime_K2}

\end{table}

As the first two sequences start with an even integer, we do not count them in $F(k)$. The third sequence, so $(1,0)$ which is generated by the integers $5 + 4q$ is such that $\chi = k = 2$. Unlike Terras, we will not count them because we have reached the condition $T^kn < n$, which will create a slight gap with the results of Terras.  Then, the distribution function $F(k)$  with $\chi > k$ instead $\chi \ge k$ really becomes the reciprocal of the function $G(k)$ defined by Lagarias. As expected, $F(k=2) = 1/4$.
 
And so on for different values of the number of iterations $k$.

The number of different sequences is given by $b = 2^k$, and the number of different parameters $a$ is calculated by the binomial coefficients $\left( \begin{array}{c} k \\ k_2 \end{array}\right)$. Binomial coefficients can be represented in a Pascal triangle,


\begin{table}[H]
\begin{center}
\begin{tabular}{c|ccccccccc}
	
	$k_2 \setminus k$  & 0 & 1 & 2 & 3 & 4 & 5  & 6  & 7  & $\cdots$ \\
		\hline
	0 & 1 & 1 & 1 & 1 & 1  & 1  & 1  & 1  & $\cdots$ \\
	1 &   & 1 & 2 & 3 & 4  & 5  & 6  & 7  & $\cdots$ \\
	2 &	  &		& 1 & 3 & 6  & 10 & 15 & 21 & $\cdots$ \\
	3	&		&		&   & 1 & 4  & 10 & 20 & 35 & $\cdots$ \\
	4	&		&		&		&   & 1  & 5  & 15 & 35 & $\cdots$ \\
	5	&		&		&		&   &    & 1  & 6  & 21 & $\cdots$ \\
	6	&		&		&		&   &    &    & 1  & 7  & $\cdots$ \\			
	7	&		&		&		&   &    &    &    & 1  & $\cdots$ \\
	8	&		&		&		&   &    &    &    &    & $\cdots$ \\ 
	

\end{tabular}
\end{center}
\caption{Pascal triangle - Binomial coefficients}
\label{PascalTriangle_BC}

\end{table}

We use a similar table which will contain the number of integers $n(i,j)$ by $2^k$ consecutive integers which satisfy the condition that the the stopping time $\chi$ is greater than the number of iterations $k$. We have 


\begin{table}[H]
\begin{center}
\begin{tabular}{c|cccccccccccc}
	
	$k_2 \setminus k$  & 0 & 1 & 2 & 3 & 4 & 5  & 6  & 7  & 8 & 9 & 10 & $\cdots$ \\
		\hline
	0   & 1 & 0 & 0 & 0 & 0  & 0  & 0  & 0  & 0  & 0  & 0  & $\cdots$ \\
	1   &   & 1 & 0 & 0 & 0  & 0  & 0  & 0  & 0  & 0  & 0  & $\cdots$ \\
	2   &	  &		& 1 & 1 & 0  & 0  & 0  & 0  & 0  & 0  & 0  & $\cdots$ \\
	3	  &		&		&   & 1 & 2  & 0  & 0  & 0  & 0  & 0  & 0  & $\cdots$ \\
	4	  &		&		&		&   & 1  & 3  & 3  & 0  & 0  & 0  & 0  & $\cdots$ \\
	5	  &		&		&		&   &    & 1  & 4  & 7  & 0  & 0  & 0  & $\cdots$ \\
	6	  &		&		&		&   &    &    & 1  & 5  & 12 & 12 & 0  & $\cdots$ \\			
	7	  &		&		&		&   &    &    &    & 1  & 6  & 18 & 30 & $\cdots$ \\
	8	  &		&		&		&   &    &    &    &    & 1  & 7  & 25 & $\cdots$ \\ 
	9	  &		&		&		&   &    &    &    &    &    & 1  & 8  & $\cdots$ \\ 
	10	&		&		&		&   &    &    &    &    &    &    & 1  & $\cdots$ \\ 		
	  	&		&		&		&   &    &    &    &    &    &    &    & $\cdots$ \\


\end{tabular}
\end{center}
\caption{Pascal triangle - Number of integers $n(i,j)$ by $2^k$ consecutive integers with $\chi > k$}
\label{PascalTriangle_nbr}

\end{table}

The index $j$ for the columns of the table is the exponent $k$ (the number of iterations) of $2$ in the parameter $b = 2^k$. The index $i$ for the rows is the exponent $k_2$ of $3$ in the parameter $a = 3^{k_2}$. As $k_2$ correspond to the number of operations on the odd integers, this value is in fact the number of $1$ in the diadic sequences and varies of $0$ to $k$. The various data in this table are calculated recursively.

The first data is trivial and indicates that all the integers satisfy the condition $\chi > k$ and this, because the number of iterations is $k = 0$. The cas $k = 1$ has ready be analyzed and we perform the following initialization, so $n(0,1) = 0$ and $n(0,1) = 1$. After $1$ iteration, all positive even integer go to a smaller integer $(n = 0)$ and, all positive odd integer go to a greater integer $(n = 1)$.   

From $k = 2$ we proceed recursively in the calculation of $n(i,k)$. 

We use the principle that each sequence is generated so that the new parameter $b$ (for $k$) is the precedent (for $k - 1$) time $2$, and the new parameter $a$ is the precedent time $1$ or $3$. 

For example, for $k = 2$, we have two $n$ which precede (for $k = 1$), so $n(0,1) = 0$ and $n(1,1) = 1$. As $n(0,1) = 0$, the sequences starting with a even positive integer for $k = 2$ will not contribute to $F(k)$ and $n(0,2) = 0$. On the other hand, the sequences generated by the integers with $n(1,1) = 1$ can contribute to $n(1,2)$ and $n(2,2)$. The new parameter $b$ will be $b = 2 \cdot 2$ and the new parameter $a$ will be $a = 3 \cdot 1$ or $a = 3 \cdot 3$ (table~\ref{stoppingTime_K2}). In the first case, $b > a$, $x \ge y$ and $\chi = k$. Then $n(1,2) = 0$. In the second case, $b < a$, $x < y$ and $\chi > k$. Then $n(2,2) = 1$. And so on for different values of $k$. 

The sum on the index $i$ of $n(i,k) / 2^k$ for a given $k$ gives the value of the distribution function $F(k)$ for this number of iterations $k$,

\begin{equation}
	F(k) = \sum_{i = 0}^{k} \frac{n(i,k)}{2^k},
\label{functionDistribution}		
\end{equation} 

It is then easy to build the computer programs starting from the recursive function worked out by Terras and by the previous process which makes it possible to fill the table~\ref{PascalTriangle_nbr}. The results of these two programs are compiled in the table~\ref{Distribution function F(k)}. We have also extended the programs to the distribution function $F_5(k)$ generated by the $5x + 1 $ function $T_5$ which we will deal with in the next section.      

Unlike to the function $F_5(k)$, the distribution function $F(k)$ decreases constantly, monotonically and asymptotically. On the other hand, $F(k)$ never become equal to zero, which is easily verified during the construction of table~\ref{PascalTriangle_nbr}. There therefore always remain positive integers which satisfy to $T^k(n) > n$, but there is an increasingly restricted although infinite set.

The process we used in the previous sections (problem $3x+1$) is somewhat different. We also use the periodicity property there, but we analyze everything in terms of geometric series and not of the evolution of a certain parameter $b$ (built from a product of 2) versus another parameter $a$ (built from a product of 3 and the 1). With groupings of integer sequences into geometric series, we can precisely determine the percentage of these which are such that $T^k(n) < n$. By continuing this process iteratively with the remaining sequences, we have shown that this percentage (always calculable) becomes smaller and smaller, and that the new remaining sequences (with $T^k(n) < n$) start with larger and larger integers. 


It should not be forgotten that these exercises have never made it possible to exclude cycles other cycles than the trivial cycle. On this last subject, we refer readers to the two previous papers~\cite{tremblay_1, tremblay_2}. 
\section{Problem 5x + 1}

Similar to the $3x + 1$ problem, we present three functions~\cite{kontorovich_lagarias} which encode the $5x + 1$ problem. 

Let the \textsl{Collatz $5x + 1$ function} $C_5(n)$ be defined as follow  

\begin{equation}
C_5(n)=\left\lbrace  
\begin{array}{ll}
5n+1 & \mbox{, if $n\equiv1\pmod{2}$}\\
\\
\frac{n}{2} & \mbox{, if $n\equiv0\pmod{2}$}.\\
\end{array}
\right.  
\label{Collatz_function_5x_plus_1}
\end{equation}

Define the \textsl{$5x + 1$ function} $T_5(n)$ (or $5x + 1$ map)

\begin{equation}
T_5(n)=\left\lbrace  
\begin{array}{ll}
\frac{5n+1}{2} & \mbox{, if $n\equiv1\pmod{2}$}\\
\\
\frac{n}{2} & \mbox{, if $n\equiv0\pmod{2}$}.\\
\end{array}
\right.  
\label{3x+1_function_5x_plus_1}
\end{equation} 

The third function, the \textsl{accelerated $5x + 1$ function} $U_5(n)$, is defined on the domain of all odd integers, and removes all powers of $2$ at each step. It is given by 

\begin{equation}
	U_5(n) = \frac{5n + 1}{2^{ord_2(5n+1)}}
\label{accelerated_function_5x_plus_1}
\end{equation}

in which $ord_2(n)$ counts the number of powers of $2$ dividing $n$. 

The function $T_5$ is known to have $5$ cycles~\cite{matthews}, with starting values $0, 1, 13, 17, -1$.

          
\subsection{Groupings (in triplets) of trajectories generated by the function $U_5(n)$}

Let us represent the \textsl{accelerated $5x + 1$ function} $U_(n)$ by the following $5$ operations

\begin{equation}
U_5(n)=\left\lbrace  
\begin{array}{lll}
\frac{5n+1}{4} & \mbox{, if $n\equiv7\pmod{8}$}\\
\\
\frac{5n+1}{2} & \mbox{, if $n\equiv1\pmod{4}$}\\
\\
\frac{5n+1}{16} & \mbox{, if $n\equiv3\pmod{32}$}\\
\\
\frac{5n+1}{8} & \mbox{, if $n\equiv11\pmod{16}$}\\
\\
\frac{5n+1}{32} & \mbox{, if $n\equiv19\pmod{32}$}\\
\\
\end{array}
\right.  
\label{accelerated_function_2_5x_plus_1}
\end{equation}

The first four operations result in odd integers and the last operation ends with a alternation of even and odd integers $3\pmod{5}$. The table~\ref{Triplet_function_U_5x_plus_1} we give the five groupings resulting of these operations including the intermediate (abreviate 'int') operation $(5x + 1 )/2$. In each grouping we have a sequence of three integers ('triplets').


          
\subsection{New function $\mathfrak{U}_5(n)$ and one-to-one correspondence with the function $U_5(n)$}

We define the function $\mathfrak{U}(n)_5$

$$\mathfrak{U}_5(n)=\left\lbrace
\begin{array}{ll}
\frac{5n}{4} & \mbox{, if $n\equiv4\pmod{4}$}\\
\\
\frac{5n-1}{2} & \mbox{, if $n\equiv1\pmod{2}$}\\
\\
\frac{5n+6}{16} & \mbox{, if $n\equiv2\pmod{16}$}\\
\\
\frac{5n+2}{8} & \mbox{, if $n\equiv6\pmod{8}$}\\
\\
\frac{n+6}{16} & \mbox{, if $n\equiv10\pmod{16}$}\\
\end{array}
\right.$$

This function does not come from generalized $3x + 1$ mappings. The table~\ref{Triplet_function_U_transformed_5x_plus_1} we give the five groupings resulting of these operations. In fact, this table brings together the trajectories of all integers and the integer resulting from a single application of the function $\mathfrak{U}_5(n)$ (duos). We form triplets by repeating the first integer. We can easily verify that all the integers $n$ of triplets in the table~\ref{Triplet_function_U_5x_plus_1} (function $U_5(n)$) are calculated using two simple operations on all the integers $n_{(new)}$, where $n_{(new)}$ is any integer beginning the triplets in the table~\ref{Triplet_function_U_transformed_5x_plus_1} (function $\mathfrak{U}_5(n)$), 

\begin{equation}
	n = 2n_{(new)} - 1 
\label{first_operation_5x_plus_1}
\end{equation} 

and 

\begin{equation}
	n = 5n_{(new)} - 2. 
\label{second_operation_5x_plus_1}
\end{equation}

All intermediate integers $3\pmod{5}$ of triplets of $U_5(n)$ and the last of the fifth grouping come from the operation $(5n_{(new)} - 2)$. All other integers (odd integers) come from the operation $(2n_{(new)} - 1)$. Then, the correspondence between the results of the application of the function $\mathfrak{U}_5(n)$ on all integers is one-to-one with the application of the function $U_5(n)$ on all odd integers.


          
\subsection{Tree generated by the function $\mathfrak{U}_5(n)$ and distribution of integers}

Let us define the branches as the sequence of the integers that are connected to each other by the operations $(5n)/4$, $(5n - 1)/2$, $(5n + 6)/16$, $(5n + 2)/8$ starting with $3\pmod{5}$ and ending with an integer resulting from the operation $(n + 6)/16$. The branches are interrelated, thus forming a tree. The $\mathfrak{U}_5(n)$ function has $5$ cycles as the function $T_5$. An infinity of integers enter these cycles, but there may be others that belong to the divergent trajectories. 

The search for the convergence towards the cycles is to determine if all the integers are in branches (distribution) interlocking into each other and if, sooner or later, they converge to smaller integers, otherwise there is divergence.

As for the function $\mathfrak{U}(n)$, all integers are found in sequences with an end of all possible lengths and their distribution is characterized by a geometric progression. 

Using the operation $(5n_{(new)} - 2)$ on all integers in branch ends of the tree and the operation $(2n_{(new)} - 1)$ on all the other integers of this tree, we build the tree generated by the function $U_5(n)$. Nevertheless, it is easier to analyze the function $\mathfrak{U}_5(n)$ which applies to all integers and which naturally reveals the regularity in the distributions of the branches. 
        

The sequences that we will analyze end with an integer resulting the operation $(n + 6)/16$, but the beginning of the sequences is any integer, not necessarily a branch beginner integer $3\pmod{5}$.   

We can write the different diophantine equations, as for the $3x + 1$ problem, and find the solutions of the latter for all the possible lengths. This approach leads us quickly to trajectories of all integers whose distribution follows a geometric progression. We will have the next results. 

The results below are valid  for $L \ge 2$.

The total number of diophantine equations $\#(DE)_L$  is

\begin{equation}
	\#(DE)_L = 4^{L-2}. 
\label{number_DE_5x_plus_1}
\end{equation}

The different increments (or modulo) for the integers $x$ starting these sequences are

\begin{equation}
	Increment = d_1*d_2*d_3*\cdots*d_{L-1},
\label{increment_5x_plus_1}
\end{equation}

where $d_1, d_2, d_3, \cdots, d_{L-2} = 2,4,8,16$ and $d_{L-1} = 16$. 

The increment (or modulo) for the last integers $y$ is $5^{L-2}$.

The different distributions of first integers $x$ are calculated from the following quadrinomial.  

\begin{equation}
	D_L = \frac{1}{16}\left(\frac{1}{2} + \frac{1}{4} + \frac{1}{8} + \frac{1}{16}\right)^{L-2},
\label{geometricSeries_firstDistribution_5x_plus_1}	
\end{equation}

which is the nth term ($n = L-1$) of a geometric progression. From equation~\ref{geometricSeries}, which is the general form for this series, we put $a = \frac{1}{16}$ and $r = \frac{15}{16}$. 

The quadrinomial coefficients $QC$ (equation~\ref{multinomial_coefficients}) of this quadrinomial give the number of different diophantine equations for a given increment.

The sum of $n$ first terms is given by the equation~\ref{geometricSum}

\begin{equation}
	S_L = \frac{1}{16}\left(\frac{(1-(\frac{15}{16})^{L-1})}{(1-\frac{15}{16})}\right) = 1-\left(\frac{15}{16}\right)^{L-1}.
\label{geometricSeries_firstDistribution_sum_5x_plus_1}	
\end{equation} 

When $L$ goes to infinity the sum tends towards one, meaning that {\bf all integers $x$, without exception, start sequences with end}. 


According to the equation~\ref{geometricSeries_firstDistribution_sum_5x_plus_1}, a large part of natural integers is included in the first sequences. The sequences with $L = 2$ until $L = 20$ include almost $70.7 \%$ of all natural integers. For this range of lengths, we covered $99.6 \%$ of all integers with the function $\mathfrak{U}(n)$. To reach this percentage with the function $\mathfrak{U}_5(n)$, it is necessary to cover the lengths $L = 2$ until $L = 85$. 


The first sequences where the final integers $y$ are larger than the first integers are those starting with $x\equiv155\pmod{256}$, $x\equiv367\pmod{512}$, $x\equiv412\pmod{512}$, $x\equiv435\pmod{512}$ and $x\equiv453\pmod{512}$, for the lengths $L = 6$. In fact, these sequences are the only ones of all those for the lengths $L = 2$ until $L = 6$ with this behavior. All other integers for sequences covering these lengths end with smaller integers. The distributions $1/256$ (or $1/512$) meaning that at each $256$  (or $512$) consecutive integers there is one that starts a sequence which end with a larger integer. 

For the sequences of length $L = 7$ and more we observe the same type of behavior. Some sequences with the smaller modulo will have the beginning integer smaller than that of the end of the sequence. For all others we have the inverse behavior.   



If we take consecutive integer slices of $5 \phantom{1} 000$, for example, we can already see that the number of integers whose last of the sequence is larger than the first is close of $60 \%$, meaning than $60 \%$ of {\bf all integers} behaving this behavior. This percentage was around $12 \%$ for the problem $3x + 1$ and the function $\mathfrak{U}(n)$. We can use the programmed function $\mathfrak{F}$ determining the fraction $f$ (or the \%) of natural integers such as $x < y$. For $L = 2$ to $20, 30, 40, 50, 60, 70, 80, 85$, we have 

\begin{flushleft}
	f = 0.3092148 (70.7 \%), 0.4455945 (84.6 \%), 0.5184832 (91.9 \%), 0.5568361 (95.8 \%),
	
	0.5769612 (97.8 \%), 0.5875168 (98.8 \%) 0,5930529 (99.4 \%),  0,5947369 (99.6 \%).  
\end{flushleft}

Between parentheses we have the proportion of all natural numbers for all sequences of length $2$ until $L$ calculated with the equation~\ref{geometricSeries_firstDistribution_sum_5x_plus_1}. Recall that we have programmed the function $\mathfrak{F}$ in VBA (Visual Basic for Applications) with a standard laptop (standard computer). From around $L = 53$ the integers appearing in the function become very large; for example, the quadrinomial coefficients have around 29 digits. By using a more powerful computer, we could obtain the results for $L > 53$.



Unlike the $3x+1$ problem, the divergence seems quite possible. For example, the trajectory generated by the function $\mathfrak{U}_5(n)$ starting with $4$ diverges quickly. Here are the first branches 

\begin{flushleft}
$4\longrightarrow5 \longrightarrow12 \longrightarrow\ldots \longrightarrow 3\phantom{1}978\phantom{1}842\longrightarrow248\phantom{1}678\quad L = 30 \quad PP$

$248\phantom{1}678 \longrightarrow 155\phantom{1}424 \longrightarrow \ldots \longrightarrow 86\phantom{1}277\phantom{1}722 \longrightarrow 5\phantom{1}392\phantom{1}358\quad L = 22 \quad PP$

$5\phantom{1}392\phantom{1}358 \longrightarrow 3\phantom{1}370\phantom{1}224 \longrightarrow \ldots \longrightarrow 957\phantom{1}875\phantom{1}242 \longrightarrow 59\phantom{1}867\phantom{1}203\quad L = 19 \quad PP$

$\ldots$

\end{flushleft}

Using the function $U_5(n)$ applied to odd integers the percentage of their distribution with $x<y$ is still higher than $60 \%$ at the first level. The value is around $68 \%$.

\section{Conclusion} 

By properly grouping the trajectories generated by the functions at the origin of the $3x+1$ and $5x+1$ problems, or rather by their accelerated functions, we could bring out the regularity that was hidden there. This is allowed us to gather the trajectories in groups of all possible lengths. Subsequently, we have been able to determine the distribution of integers $x$ that start these trajectories and, finally, to find the proportions of them that are smaller (or larger) than the integer $y$ ending the trajectory. The distribution of integers following a geometric progression has given us the opportunity to treat all integers, without exception, until infinity. In no case did we use the notions of probabilities (or stochastic process), which are widely used in the literature. The path we followed is clearly different from all that has been done nowadays to resolve the conjectures and brings a whole new light. We believe that a new way is open for the comprehension and the solution of the behavior of many other functions that are parted of the generalized $3x+1$ mappings, or others.

We conclude this exercise by quoting a sentence taken from the conclusion from a $1985$ paper written by Lagarias~\cite{lagarias}, namely :  "Of course there remains the possibility that someone will find some hidden regularity in the $3x+1$ problem that allows some of the conjectures about it to be settled". The theorem on the periodicity bring out the regularity in the trajectories generated by the $3x + 1$ function $T(n)$, the $5x + 1$ function $T_5(n)$ and many other functions. The application of this theorem to the various problems generated by the functions derived from generalized $3x + 1$ mappings is surely an essential key to the implementation of their solutions together with the basic properties of the diophantine equations.    




\newpage

\begin{table}
\begin{center}
\begin{tabular}{*{17}{|l}|}
	
	\hline
		\multicolumn{5}{|c|}{\bf $x_0\equiv1\pmod{8}$} & & \multicolumn{5}{c|}{\bf $x_0\equiv3\pmod{4}$} & &
		\multicolumn{5}{c|}{\bf $x_0\equiv5\pmod{8}$} \\
	\hline
		\multicolumn{5}{|c|}{$\hdots$} & & \multicolumn{5}{c|}{$\hdots$} & &
		\multicolumn{5}{c|}{$\hdots$} \\
	\hline
		1 & $\rightarrow$ & \cellcolor[gray]{0.8} 2 & $\rightarrow$ & 1 & & 3 & $\rightarrow$ & \cellcolor[gray]{0.8} 5 & $\rightarrow$ & 5 &  & 5  & $\rightarrow$ & \cellcolor[gray]{0.8} 8 & $\rightarrow$ & 2 \\
	\hline
		9 & $\rightarrow$ & \cellcolor[gray]{0.8} 14 & $\rightarrow$ & 7 & & 7 & $\rightarrow$ & \cellcolor[gray]{0.8} 11 & $\rightarrow$ & 11 &  & 13  & $\rightarrow$ & \cellcolor[gray]{0.8} 20 & $\rightarrow$ & 5 \\
	\hline
		17 & $\rightarrow$ & \cellcolor[gray]{0.8} 26 & $\rightarrow$ & 13 & & 11 & $\rightarrow$ & \cellcolor[gray]{0.8} 17 & $\rightarrow$ & 17 &  & 21  & $\rightarrow$ & \cellcolor[gray]{0.8} 32 & $\rightarrow$ & 8 \\
	\hline
		25 & $\rightarrow$ & \cellcolor[gray]{0.8} 38 & $\rightarrow$ & 19 & & 15 & $\rightarrow$ & \cellcolor[gray]{0.8} 23 & $\rightarrow$ & 23 &  & 29  & $\rightarrow$ & \cellcolor[gray]{0.8} 44 & $\rightarrow$ & 11 \\
	\hline
		33 & $\rightarrow$ & \cellcolor[gray]{0.8} 50 & $\rightarrow$ & 25 & & 19 & $\rightarrow$ & \cellcolor[gray]{0.8} 29 & $\rightarrow$ & 29 &  & 37  & $\rightarrow$ & \cellcolor[gray]{0.8} 56 & $\rightarrow$ & 14 \\
	\hline
		41 & $\rightarrow$ & \cellcolor[gray]{0.8} 62 & $\rightarrow$ & 31 & & 23 & $\rightarrow$ & \cellcolor[gray]{0.8} 35 & $\rightarrow$ & 35 &  & 45  & $\rightarrow$ & \cellcolor[gray]{0.8} 68 & $\rightarrow$ & 17 \\
\hline
	
	\end{tabular}
\end{center}
\caption{Groupings in triplets of trajectories generated by the function $U(n)$}
\label{Triplet_function_U}
\end{table}

\begin{table}
\begin{center}
\begin{tabular}{*{17}{|l}|}
	
	\hline
		\multicolumn{5}{|c|}{\bf ${x_0}_{(new)}\equiv1\pmod{4}$} & & \multicolumn{5}{c|}{\bf ${x_0}_{(new)}\equiv2\pmod{2}$} & &
		\multicolumn{5}{c|}{\bf ${x_0}_{(new)}\equiv3\pmod{4}$} \\
	\hline
		\multicolumn{5}{|c|}{$\hdots$} & & \multicolumn{5}{c|}{$\hdots$} & &
		\multicolumn{5}{c|}{$\hdots$} \\
	\hline
		1 & $\rightarrow$ & \cellcolor[gray]{0.8} 1 & $\rightarrow$ & 1 & & 2 & $\rightarrow$ & \cellcolor[gray]{0.8} 2 & $\rightarrow$ & 3 &  & 3  & $\rightarrow$ & \cellcolor[gray]{0.8} 3 & $\rightarrow$ & 1 \\
	\hline
		5 & $\rightarrow$ & \cellcolor[gray]{0.8} 5 & $\rightarrow$ & 4 & & 4 & $\rightarrow$ & \cellcolor[gray]{0.8} 4 & $\rightarrow$ & 6 &  & 7  & $\rightarrow$ & \cellcolor[gray]{0.8} 7 & $\rightarrow$ & 2 \\
	\hline
		9 & $\rightarrow$ & \cellcolor[gray]{0.8} 9 & $\rightarrow$ & 7 & & 6 & $\rightarrow$ & \cellcolor[gray]{0.8} 6 & $\rightarrow$ & 9 &  & 11  & $\rightarrow$ & \cellcolor[gray]{0.8} 11 & $\rightarrow$ & 3 \\
	\hline
		13 & $\rightarrow$ & \cellcolor[gray]{0.8} 13 & $\rightarrow$ & 10 & & 8 & $\rightarrow$ & \cellcolor[gray]{0.8} 8 & $\rightarrow$ & 12 &  & 15  & $\rightarrow$ & \cellcolor[gray]{0.8} 15 & $\rightarrow$ & 4 \\
	\hline
		17 & $\rightarrow$ & \cellcolor[gray]{0.8} 17 & $\rightarrow$ & 13 & & 10 & $\rightarrow$ & \cellcolor[gray]{0.8} 10 & $\rightarrow$ & 15 &  & 19  & $\rightarrow$ & \cellcolor[gray]{0.8} 19 & $\rightarrow$ & 5 \\
	\hline
		21 & $\rightarrow$ & \cellcolor[gray]{0.8} 21 & $\rightarrow$ & 16 & & 12 & $\rightarrow$ & \cellcolor[gray]{0.8} 12 & $\rightarrow$ & 18 &  & 23  & $\rightarrow$ & \cellcolor[gray]{0.8} 23 & $\rightarrow$ & 6 \\
\hline
	
	\end{tabular}
\end{center}
\caption{Groupings in triplets of trajectories generated by the function $\mathfrak{U}$}
\label{Triplet_function_U_transformed}
\end{table}

\begin{table}
\begin{center}
\scriptsize
\rotatebox{90}{
\begin{tabular}{*{21}{p{0.42cm}}}

	&& & \cellcolor[gray]{0.8} $\dashrightarrow$ & \cellcolor[gray]{0.8} 80 &&&&&&& \cellcolor[gray]{0.8} 14 &  \cellcolor[gray]{0.8} $\dashleftarrow$ &&&&&&
	&& \\
	&& &  & \cellcolor[gray]{0.8} $\downarrow$ &&&&&&& \cellcolor[gray]{0.8} $\cdots$ &  &  &&&&&&&\\
	&& & \cellcolor[gray]{0.8} $\dashrightarrow$ & \cellcolor[gray]{0.8} 120 & &&&&&& \cellcolor[gray]{0.8} $\downarrow$ &  & &&&&&&&  \\
	&& &  & \cellcolor[gray]{0.8} $\downarrow$ &&& \cellcolor[gray]{0.8} 56 & \cellcolor[gray]{0.8} $\leftarrow$ & \cellcolor[gray]{0.8} $\leftarrow$  & \cellcolor[gray]{0.8} $\leftarrow$  & \cellcolor[gray]{0.8} 56 &  &  &&&&&&& \\
	&& & \cellcolor[gray]{0.8} $\dashrightarrow$ & \cellcolor[gray]{0.8} 180 &&& \cellcolor[gray]{0.8} $\cdots$ &&&&  &  & &&&&&&&   \\
	&& & & \cellcolor[gray]{0.8} $\downarrow$ &&& \cellcolor[gray]{0.8} $\downarrow$ &&&&  &  & &&&&&&&   \\
	&& & \cellcolor[gray]{0.8} $\dashrightarrow$ & \cellcolor[gray]{0.8} 270 &&& \cellcolor[gray]{0.8} 810 & \cellcolor[gray]{0.8} $\dashleftarrow$ &  & &  &  & & \cellcolor[gray]{0.8} $\dashrightarrow$ & \cellcolor[gray]{0.8} 641 &&&&  &\\
	&& &  & \cellcolor[gray]{0.8} $\downarrow$ &&& \cellcolor[gray]{0.8} $\downarrow$ &  &  &&  &  &&& \cellcolor[gray]{0.8} $\downarrow$ &&&&  &  \\
\hline			
\end{tabular}
}
\normalsize
\end{center}
\caption{Tree - function $\mathfrak{U}$}
\label{Tree}
\end{table}

\begin{table}
\begin{center}
\scriptsize
\rotatebox{90}{
\begin{tabular}{*{22}{p{0.405cm}}}

	&& &  & \cellcolor[gray]{0.8} 159 &&&&&&& \cellcolor[gray]{0.8} 27 &   &&&&&&
	&& \\
	&& &  & \cellcolor[gray]{0.8} $\downarrow$ &&&&&&& \cellcolor[gray]{0.8} $\cdots$ & \cellcolor[gray]{0.8} $\dashleftarrow$ &  &&&&&&&\\
	&& & \cellcolor[gray]{0.8} $\dashrightarrow$ & \cellcolor[gray]{0.8} 239 & &&&&&& \cellcolor[gray]{0.8} $\downarrow$ &  & &&&&&&&  \\
	&& &  & \cellcolor[gray]{0.8} $\downarrow$ &&& \cellcolor[gray]{0.8} 111 &  &   &   & \cellcolor[gray]{0.8} $\downarrow$ &  &  &&&&&&& \\
	&& & \cellcolor[gray]{0.8} $\dashrightarrow$ & \cellcolor[gray]{0.8} 359 &&& \cellcolor[gray]{0.8} $\cdots$ & \cellcolor[gray]{0.8} $\leftarrow$ & \cellcolor[gray]{0.8} $\leftarrow$ & \cellcolor[gray]{0.8} $\leftarrow$ & \cellcolor[gray]{0.8} 167 &  & &&&&&&&   \\
	&& & & \cellcolor[gray]{0.8} $\downarrow$ &&& \cellcolor[gray]{0.8} $\downarrow$ &&&&  &  & &&&&&&&   \\
	&& & \cellcolor[gray]{0.8} $\dashrightarrow$ & \cellcolor[gray]{0.8} 539 &&& \cellcolor[gray]{0.8} 1619 &  &  & &  &  & &  & \cellcolor[gray]{0.8} 1281 &&&&  &\\
	&& &  & \cellcolor[gray]{0.8} $\downarrow$ &&& \cellcolor[gray]{0.8} $\downarrow$ &  &  &&  &  && \cellcolor[gray]{0.8} $\dashrightarrow$ & \cellcolor[gray]{0.8} $\downarrow$ &&&&  &  \\
\hline			
\end{tabular}
}
\normalsize
\end{center}
\caption{Tree - function U - odd integers}
\label{Tree_odds}
\end{table}

\begin{table}
\begin{center}
\begin{tabular}{*{1}{l}}
	L = 2 \quad 3$\rightarrow$1 \quad modulo= 4 PG\\

	\\

	L = 3 \quad 2$\rightarrow$3$\rightarrow$1 \quad modulo= 8 PG\\
	L = 3 \quad 9$\rightarrow$7$\rightarrow$2 \quad modulo= 16 PG\\

	\\

	L = 4 \quad 12$\rightarrow$18$\rightarrow$27$\rightarrow$7 \quad modulo= 16 PG\\
	L = 4 \quad 6$\rightarrow$9$\rightarrow$7$\rightarrow$2 \quad modulo= 32 PG\\
	L = 4 \quad 13$\rightarrow$10$\rightarrow$15$\rightarrow$4 \quad modulo= 32 PG\\
\end{tabular}
\end{center}
\caption{Sequences generated by the function $\mathfrak{U}$ for $L = 2$ until $L = 6$}
\label{Sequences_L}
\end{table}

\begin{table}
\begin{center}
\begin{tabular}{*{1}{l}}
	L = 2 \quad 5$\rightarrow$2 \quad modulo= 8 PG\\

	\\

	L = 3 \quad 3$\rightarrow$5$\rightarrow$2 \quad modulo= 16 PG\\
	L = 3 \quad 17$\rightarrow$13$\rightarrow$5 \quad modulo= 32 PG\\

	\\

	L = 4 \quad 23$\rightarrow$35$\rightarrow$53$\rightarrow$20 \quad modulo= 32 PG\\
	L = 4 \quad 11$\rightarrow$17$\rightarrow$13$\rightarrow$5 \quad modulo= 64 PG\\
	L = 4 \quad 25$\rightarrow$19$\rightarrow$29$\rightarrow$11 \quad modulo= 64 PG\\
\end{tabular}
\end{center}
\caption{Sequences generated by the function $U(n)$ for $L = 2$ until $L = 6$}
\label{Sequences_L_odds}
\end{table}

\begin{table}
\begin{center}

\begin{tabular}{*{8}{|l}|}
	\hline
	 $L$ & $\#(DE)$ &   \multicolumn {2}{p{1.6cm}|}{Modulo} & Dist.  & $D_L$ & $S_L$ & $S_L$\\
	  & $= 2^{L-2}$ &   \multicolumn {2}{p{1.6cm}|}{} & of integers & $=(1/4)(3/4)^{L-2}$ & $=1-(3/4)^{L-1}$ & $(\%)$ \\
	\hline
		&   &   & $BC_{n,k}$  &   &   &  & \\
	\hline	
		2 &  1 &  mod 4  & 1  & 1/4  & 1/4  & 1/4 & 25\\
	\hline	
		3 &  2 &  mod 8 & 1  & 1/8  &  & & \\
		  &    &   mod 16 & 1  & 1/16  & 3/16 & 7/16 & 43.75 \\
	\hline
		4 &  4 &  mod 16 & 1  & 1/16  &  & & \\
		 &   &  mod 32 & 2  & 2/32  &  & & \\
		 &   &  mod 64 & 1  & 1/64  & 9/64 & 37/64 & 57.81 \\
	\hline
		5 &  8 &  mod 32 & 1  & 1/32  &  & & \\
		 &   &  mod 64 & 3  & 3/64  &  & & \\
		 &   &  mod 128 & 3  & 3/128  &  & & \\
		 &   &  mod 256 & 1  & 1/256  & 27/256 & 175/256 & 68.36 \\
	\hline
		6 &  16 &  mod 64 & 1  & 1/64  &  & & \\
		 &   &  mod 128 & 4  & 4/128  &  & & \\
		 &   &  mod 256 & 6  & 6/256  &  & & \\
		 &   &  mod 512 & 4  & 4/512  &  & & \\
		 &   &   mod 1024 & 1  & 1/1024  & 81/1024 & 781/1024 & 76.27 \\
	\hline
		7 &  32 &  mod 128 & 1  & 1/128  &  & & \\
		 &   &  mod 256 & 5  & 5/256  &  & & \\
		 &   &  mod 512 & 10  & 10/512  &  & & \\
		 &   &  mod 1024 & 10  & 10/1024  &  & & \\
		 &   &   mod 2048 & 5  & 5/2048  &  &  &  \\
		 &   &   mod 4096 & 1  & 1/4096  & 243/4096 & 3367/4096 & 82.20 \\
	\hline
		8 &  64 &  $\ldots$ & $\ldots$  & $\ldots$  & $\ldots$ & $\ldots$ & $\ldots$ \\
	\hline
		$\ldots$ &  $\ldots$ &  $\ldots$ & $\ldots$  & $\ldots$  & $\ldots$ & $\ldots$ & $\ldots$ \\
	\hline		

	\end{tabular}
	
	\end{center}
	\caption{Distribution of first integers $x$ for $L = 2$ until $L = 7$ (function $\mathfrak{U}$)}
	\label{Distribution}
	\end{table}

\begin{table}
\begin{center}

\begin{tabular}{*{6}{|l}|}
	\hline
	 $L$  & $k_3 = L - 2$ & $(k_1, k_2)$ &  $b = 4 \cdot 2^{k_1} \cdot 4^{k_2}$ & $a = 3^{k_3} = 3^{L -2}$ & $x$ versus $y$ \\

	\hline	
		2 & 0 & $(0, 0)$ & 4 &  1 & $x > y$\\
	\hline	
		3 & 1 & $(1, 0)$ &  8 & 3 & $x > y$ \\
		  &  & $(0, 1)$ &16  &  & $x > y$ \\
	\hline
		4  & 2 & $(2, 0)$ & 16 &  $3^2 = 9$ & $x > y$ \\
		 & & $(1, 1)$ & 32 &  & $x > y$ \\
		 & & $(0, 2)$ & 64 &  & $x > y$\\
	\hline
		5  & 3 & $(3, 0)$ & 32 &  $3^3 = 27$ & $x > y$ \\
		 & & $(2, 1)$ & 64 &  & $x > y$ \\
		 & & $(1, 2)$ & 128 &  & $x > y$ \\
		 & & $(0, 3)$ & 256  &  & $x > y$ \\
	\hline
		6 & 4 & $(4, 0)$ & 64 &  $3^4 = 81$ & \cellcolor[gray]{0.8}$x < y$ \\
		 & & $(3, 1)$ & 128 &  & $x > y$ \\
		 & & $(2, 2)$ & 256 &  & $x > y$ \\
		 & & $(1, 3)$ & 512 &  & $x > y$ \\
		 & & $(0, 4)$ &1024 &  & $x > y$ \\
	\hline
		7  & 5 & $(5, 0)$ & 128 &  $3^5 =243$ & \cellcolor[gray]{0.8}$x < y$ \\
		 & & $(4, 1)$ & 256 &  & $x > y$ \\
		 & & $(3, 2)$ & 512 &  & $x > y$ \\
		 & & $(2, 3)$ & 1024 &  & $x > y$ \\
		 & & $(1, 4)$ & 2048 &  & $x > y$  \\
		 & & $(0, 5)$ & 4096 &  & $x > y$ \\
	\hline
		8  & 6 & $(6, 0)$ & 256 &  $3^6 = 727$ & \cellcolor[gray]{0.8}$x < y$   \\
		 & & $(5, 1)$ & 512 &  & \cellcolor[gray]{0.8}$x < y$ \\
		 & & $(4, 2)$ & 1024 &  & $x > y$ \\
		 & & $(3, 3)$ & 2048 &  & $x > y$ \\
		 & & $(2, 4)$ & 4096 &  & $x > y$ \\
		 & & $(1, 5)$ & 8192 &  & $x > y$ \\
		 & & $(0, 6)$ & 16384 &  & $x > y$ \\
	\hline
	  9  & 7 & $(7, 0)$ & 512 &  $3^7 = 2181$ & \cellcolor[gray]{0.8}$x < y$   \\
		$\ldots$ &  $\ldots$ &  $\ldots$ & $\ldots$ & $\ldots$ & $\ldots$   \\
	\hline		

	\end{tabular}
	
	\end{center}
	\caption{Distribution of parameters $b$ and $a$ for $L = 2$ until $L = 8$ (function $\mathfrak{U}$)}
	\label{Distribution_x_y}
	\end{table}
	
\begin{table}
\begin{center}
\tiny

\begin{tabular}{*{9}{|l}|}
	\hline
	 $L$ & $\#(DE)$ &   \multicolumn {2}{p{1.6cm}|}{Modulo} & Dist.  & $D_{x<y}$ & Cumul & Cumul & $S_L$ \\
	  & $= 2^{L-2}$ &   \multicolumn {2}{p{1.6cm}|}{} & of integers &  & & $(\%)$ & $(\%)$ \\
	\hline
		&  &   &  $BC_{n,k}$ &   &   &  & & \\
	\hline	
		6 &  16 &  mod 64  & 1  & 1/64  & 1/64  & 1/64 & 1.56 & 76.27 \\
	\hline	
		7 &  32 &  mod 128 & 1  & 1/128  & 1/128 & 3/128 & 2.34 & 82.20 \\		
	\hline
		8 &  64 &  mod 256 & 1  & 1/256  &  & & & \\
		 &   &  mod 512 & 6  & 6/512  & 8/512 & 20/512 & 3.91 & 86.65 \\		 
	\hline
		9 &  128 &  mod 512 & 1  & 1/512  &  & & & \\
\hline
	
	\end{tabular}
	
	\normalsize
	\end{center}
	\caption{Distribution of integers for $x < y$ with $L = 6$ until $L = 20$ (function $\mathfrak{U}$)}
	\label{Distribution_smaller}
	\end{table}

\begin{table}
\begin{center}
\begin{tabular}{*{1}{l}}

61$\rightarrow$46$\rightarrow$69$\rightarrow$52$\rightarrow$78$\rightarrow$117$\rightarrow$88$\rightarrow$132$\rightarrow$198$\rightarrow$297$\rightarrow$223$\rightarrow$56\quad L = 12\quad PG \\

142$\rightarrow$213$\rightarrow$160$\rightarrow$240$\rightarrow$360$\rightarrow$540$\rightarrow$810$\rightarrow$1215$\rightarrow$304\quad L = 9\quad PP \\

223$\rightarrow$56\quad L = 2\quad PG \\

304$\rightarrow$456$\rightarrow$684$\rightarrow$1026$\rightarrow$1539$\rightarrow$385\quad L = 6\quad PP \\

385$\rightarrow$289$\rightarrow$217$\rightarrow$163$\rightarrow$41\quad L = 5\quad PG \\

466$\rightarrow$699$\rightarrow$175\quad L = 3\quad PG \\

547$\rightarrow$137\quad L = 2\quad PG \\

628$\rightarrow$942$\rightarrow$1413$\rightarrow$1060$\rightarrow$1590$\rightarrow$2385$\rightarrow$1789$\rightarrow$1342$\rightarrow$2013$\rightarrow$1510$\rightarrow$2265$\rightarrow$1699$\rightarrow$425\quad L = 13\quad PG \\
\end{tabular}
\end{center}
\caption{Second level sequences $x_2\equiv61\pmod{81}$ (function $\mathfrak{U}$)}
\label{Second_level_sequences}
\end{table}

\newpage

\begin{table}
\begin{center}

\begin{tabular}{*{6}{|l}|}

	\hline
		level & \multicolumn{4}{|c|}{integers beginning with 1,2,3, $\cdots$} & f generated by the   \\
		& \multicolumn{4}{|c|}{} & programmed function $\mathfrak{F}$  \\	
		& \multicolumn{4}{|c|}{} & (equations \ref{fraction_level_1},\ref{fraction_level_2},\ref{fraction_level_3})  \\			
	\hline	
	\hline	
		& until 5 000 & until 10 000 & until 100 000 & until 1 000 000 & L1,L2,L3 = 2 to 30\\
	\hline
		& & & & & \\
	\hline
		1	& 0.123049 & 0.124325 & 0.123972 & 0.123814 & 0.1236245 \\
	\hline
		2	& 0.053421 & 0.051510 & 0.050701 & 0.051026 & 0.0508968\\		
	\hline
		3	& 0.025010 & 0.024105 & 0.024040 & 0.024105 & 0.0240408\\
	\hline
		4	& 0.016407 & 0.015603 & 0.012790 & 0.012441 & $\cdots$\\		
	\hline
		5	& 0.007203 & 0.007201 & 0.007220 & 0.006766 & $\cdots$\\	
	\hline
		6	& 0.003201 & 0.002801 & 0.004010 & 0.003908 & $\cdots$\\
	\hline
		7	& 0.003001 & 0.002100 & 0.002390 & 0.002389 & $\cdots$\\	
				
	\hline
		$\cdots$ & $\cdots$ & $\cdots$ & $\cdots$ & $\cdots$ & $\cdots$ \\
	\hline				
	\end{tabular}
	
	\end{center}
	\caption{Distribution of integers with $x < y$ (function $\mathfrak{U}$) by slices and anticipated by the programmed function $\mathfrak{F}$} 
	\label{Distribution_by_slices}
	\end{table}

\begin{table}
\begin{center}

\begin{tabular}{*{6}{|l}|}

	\hline
		level & \multicolumn{4}{|c|}{odd integers beginning with 1,3,5, $\cdots$} & f generated by the   \\
	
		& \multicolumn{4}{|c|}{} & programmed function $\mathfrak{F}$  \\	
	\hline
	\hline
		& until 5 001 &  until 10 001 & until 100 001 & until 1 000 001 & L1,L2,L3 = 2 to 30\\
	\hline
		& & & & & \\
	\hline
		1	& 0.200480 & 0.201640 & 0.202664 & 0.202652 & 0.2024107 \\
	\hline
		2	& 0.089236 & 0.085617 & 0.084002 & 0.084262 & 0.0840688  \\		
	\hline
		3	& 0.047619 & 0.041208 & 0.038301 & 0.039736 & 0.0397817  \\
	\hline
		4	& 0.029212 & 0.025205 & 0.020580 & 0.020810 & $\cdots$ \\		
	\hline
		5	& 0.014806 & 0.012603 & 0.011900 & 0.011700 & $\cdots$ \\	
	\hline
		6	& 0.006403 & 0.004409 & 0.006500 & 0.006608 & $\cdots$ \\
	\hline
		7	& 0.006002 & 0.004201 & 0.003780 & 0.004088 & $\cdots$ \\
				
	\hline
		$\cdots$ & $\cdots$ & $\cdots$ & $\cdots$ & $\cdots$ & $\cdots$ \\
	\hline				
	\end{tabular}
	
	\end{center}
	\caption{Distribution of odd integers with $x < y$ (function $U$) by slices and anticipated by the programmed function $\mathfrak{F}$}
	\label{Distribution_by_slices_U}
	\end{table}

\begin{table}
\begin{center}
\large


\begin{tabular}{*{5}{|l}|}
	\hline
	 $k$ &   \multicolumn {2}{c|}{F(k)} & \multicolumn {2}{c|}{$F_5(k)$} \\
	\hline
		& $Terras$ & $revisited$ & $Terras$ & $revisited$\\
		&  & $(with \phantom{1} the \phantom{1} table~\ref{PascalTriangle_nbr})$ &  & $(as \phantom{1} the \phantom{1} table~\ref{PascalTriangle_nbr})$\\
	\hline	

		0 &  1 &  1 & 1 & 1 \\
	\hline	
		1 &  0.5 &  0.5 & 0.5 & 0.5 \\	
	\hline
		2 &  0.5 &  0.25 & 0.5 & 0.5   \\ 
	\hline
		3 &  0.25 &  0.25 & 0.375 & 0.5  \\
	\hline
		4 &  0.25 &  0.1875 & 0.375 & 0.375 \\
	\hline
		5 &  0.1875 &  0.125 & 0.3125 & 0.375 \\
	\hline
		6 &  0.125 &  0.125 & 0.3125 & 0.3125 \\
	\hline		
		7 &  0.125 &  0.1015625 & 0.2734375 & 0.3125 \\
	\hline	
		8 &  0.1015625 &  0.07421875 & 0.2734375 & 0.2734375 \\
	\hline			
		9 &  0.07421875 &  0.07421875 & 0.2734375 & 0.2734375 \\
	\hline
		10 &  0.07421875 &  0.0625 & 0.259765625 & 0.2734375 \\
	\hline
		20 &  $2.8591 \times 10^{-2}$ &  $2.6062 \times 10^{-2}$ & 0.2221222 & 0.2221222 \\
	\hline
		30 &  $1.1894 \times 10^{-2}$ &  $1.1894 \times 10^{-2}$ & 0.2057265 & 0.2057265 \\
	\hline	
		40 &  $6.5693 \times 10^{-3}$ &  $5.8233 \times 10^{-3}$ & 0.19625785 & 0.1978474 \\
	\hline		
		50 &  $3.5373 \times 10^{-3}$ &  $3.3167 \times 10^{-3}$ & 0.191165633 & 0.191165633 \\
	\hline
		60 &  $1.9222 \times 10^{-3}$ &  $1.9222 \times 10^{-3}$ & 0.18811449 & 0.18811449 \\
	\hline	
		70 &  $1.1644 \times 10^{-3}$ &  $1.0516 \times 10^{-3}$ & 0.18513014 & 0.185735 \\
	\hline	
		80 &  $7.0744 \times 10^{-4}$ &  $6.6440 \times 10^{-4}$ & 0.1831778 & 0.1831778 \\
	\hline	
		90 &  $4.1078 \times 10^{-4}$ &  $4.1078 \times 10^{-4}$ & 0.1819218 & 0.1819218 \\
	\hline	
		100 &  $2.6396 \times 10^{-4}$ &  $2.3868 \times 10^{-4}$ & & \\
	\hline	
	
	  $\ldots$ &  $\ldots$ &  $\ldots$ &  $\ldots$ &  $\ldots$ \\ 
	\hline
	
	\end{tabular}
	
	\normalsize
	\end{center}
	\caption{Distribution function F(k)}
	\label{Distribution function F(k)}
	\end{table}

\begin{table}
\begin{center}
\begin{tabular}{*{15}{|l}|}
	
	\hline
		\multicolumn{3}{|c|}{\bf $x_0\equiv7\pmod{8}$}& \multicolumn{3}{c|}{\bf $x_0\equiv1\pmod{4}$}&
		\multicolumn{3}{c|}{\bf $x_0\equiv3\pmod{32}$} & \multicolumn{3}{c|}{\bf $x_0\equiv11\pmod{16}$}&
		\multicolumn{3}{c|}{\bf $x_0\equiv19\pmod{32}$} \\
	\hline
		\multicolumn{3}{|c|}{$\hdots$} & \multicolumn{3}{c|}{$\hdots$} &
		\multicolumn{3}{c|}{$\hdots$} & \multicolumn{3}{c|}{$\hdots$} & \multicolumn{3}{c|}{$\hdots$} \\
	\hline
		7 & \cellcolor[gray]{0.8} 18 & 9 & 1 & \cellcolor[gray]{0.8} 3 & 3 & 3 & \cellcolor[gray]{0.8} 8 & 1 & 11 & \cellcolor[gray]{0.8} 28 & 7 & 19 & \cellcolor[gray]{0.8} 48 & 3\\
	\hline
		15 & \cellcolor[gray]{0.8} 38 & 19 & 5 & \cellcolor[gray]{0.8} 13 & 13 & 35 & \cellcolor[gray]{0.8} 88 & 11 & 27 & \cellcolor[gray]{0.8} 68 & 17 & 51 & \cellcolor[gray]{0.8} 128 & 8\\
	\hline
		23 & \cellcolor[gray]{0.8} 58 & 29 & 9 & \cellcolor[gray]{0.8} 23 & 23 & 67 & \cellcolor[gray]{0.8} 168 & 21 & 43 & \cellcolor[gray]{0.8} 108 & 27 & 83 & \cellcolor[gray]{0.8} 208 & 13\\
	\hline
		31 & \cellcolor[gray]{0.8} 78 & 39 & 13 & \cellcolor[gray]{0.8} 33 & 33 & 99 & \cellcolor[gray]{0.8} 248 & 31 & 59 & \cellcolor[gray]{0.8} 148 & 37 & 115 & \cellcolor[gray]{0.8} 288 & 18\\
	\hline
		39 & \cellcolor[gray]{0.8} 98 & 49 & 17 & \cellcolor[gray]{0.8} 43 & 43 & 131  & \cellcolor[gray]{0.8} 328 & 41 & 75 & \cellcolor[gray]{0.8} 188 & 47 & 147 & \cellcolor[gray]{0.8} 368 & 23\\
	\hline
		47 & \cellcolor[gray]{0.8} 118 & 59 & 21 & \cellcolor[gray]{0.8} 58 & 53 & 163 & \cellcolor[gray]{0.8} 408 & 51 & 91 & \cellcolor[gray]{0.8} 228 & 57 & 179 & \cellcolor[gray]{0.8} 448 & 28\\
\hline
	
	\end{tabular}
\end{center}
\caption{Groupings in triplets of trajectories generated by the function $U_5(n)$}
\label{Triplet_function_U_5x_plus_1}
\end{table}

\begin{table}
\begin{center}
\begin{tabular}{*{15}{|l}|}
	
	\hline
		\multicolumn{3}{|c|}{\bf ${x_0}\equiv4\pmod{4}$} & \multicolumn{3}{c|}{\bf ${x_0}\equiv1\pmod{2}$} &
		\multicolumn{3}{c|}{\bf ${x_0}\equiv2\pmod{16}$} & \multicolumn{3}{c|}{\bf ${x_0}\equiv6\pmod{8}$} &
		\multicolumn{3}{c|}{\bf ${x_0}\equiv10\pmod{16}$}\\
	\hline
		\multicolumn{3}{|c|}{$\hdots$} & \multicolumn{3}{c|}{$\hdots$} &
		\multicolumn{3}{c|}{$\hdots$} & \multicolumn{3}{c|}{$\hdots$} &
		\multicolumn{3}{c|}{$\hdots$}\\
	\hline
		4 & \cellcolor[gray]{0.8} 4 & 5 & 1 & \cellcolor[gray]{0.8} 1 & 2 & 2 & \cellcolor[gray]{0.8} 2 & 1 & 6 & \cellcolor[gray]{0.8} 6 & 4 &  10  & \cellcolor[gray]{0.8} 10  & 1\\
	\hline
		8 & \cellcolor[gray]{0.8} 8 & 10 & 3 & \cellcolor[gray]{0.8} 3 & 7 & 18 & \cellcolor[gray]{0.8} 18 & 6 & 14 & \cellcolor[gray]{0.8} 14  & 9 & 26  & \cellcolor[gray]{0.8} 26 & 2\\
	\hline
		12 & \cellcolor[gray]{0.8} 12 & 15 & 5 & \cellcolor[gray]{0.8} 5 & 12 & 34 & \cellcolor[gray]{0.8} 34 & 11 & 22 & \cellcolor[gray]{0.8} 22 & 14 & 42 & \cellcolor[gray]{0.8} 42 & 3\\
	\hline
		16  & \cellcolor[gray]{0.8} 16 & 20 & 7 & \cellcolor[gray]{0.8} 7 & 17 & 50 & \cellcolor[gray]{0.8} 50 & 16 & 30 & \cellcolor[gray]{0.8} 30 & 19 & 58 & \cellcolor[gray]{0.8} 58 & 4\\
	\hline
		20 & \cellcolor[gray]{0.8} 20 & 25 & 9 & \cellcolor[gray]{0.8} 9 & 22 & 66 & \cellcolor[gray]{0.8} 66 & 21 & 38 & \cellcolor[gray]{0.8} 38 & 24 & 74 & \cellcolor[gray]{0.8} 74 & 5\\
	\hline
		24 & \cellcolor[gray]{0.8} 24 & 30 & 11 & \cellcolor[gray]{0.8} 11 & 27 & 82 & \cellcolor[gray]{0.8} 82 & 66 & 46 & \cellcolor[gray]{0.8} 46 & 29 & 90 & \cellcolor[gray]{0.8} 90 & 6\\
\hline
	
	\end{tabular}
\end{center}
\caption{Groupings in triplets of trajectories generated by the function $\mathfrak{U}_5$}
\label{Triplet_function_U_transformed_5x_plus_1}
\end{table}

\end{document}